\documentclass[a4paper,leqno]{article}

\usepackage[english]{babel}
\usepackage[utf8x]{inputenc}
\usepackage{enumerate,amsmath,amssymb,amstext,amsthm}
\usepackage{graphicx}
\usepackage[colorinlistoftodos]{todonotes}
\usepackage{authblk}

\usepackage{hyperref}

\newtheorem{theorem}{Theorem}[section]
\newtheorem{lemma}[theorem]{Lemma}
\newtheorem{corollary}[theorem]{Corollary}
\newtheorem{proposition}[theorem]{Proposition}

\theoremstyle{definition}

\newtheorem{defn}[theorem]{Definition}
\newtheorem{remark}[theorem]{Remark}

\def\packed{\setlength{\itemsep}{1pt} \setlength{\parskip}{0pt} \setlength{\parsep}{0pt} }

\newcommand{\R}{\mathbb{R}}

\DeclareMathOperator{\spn}{span}

\newcommand{\one}{\ensuremath{\vec{1}}}

\newcommand{\chiv}{\chi_{v}}
\newcommand{\chisv}{\chi_{sv}}

\newcommand{\transpose}{^T}
\newcommand{\spec}{\mathcal{S}(G,\p)}

\newcommand{\p}{\text{{\bf{p}}}}
\newcommand{\q}{\text{{\bf{q}}}}

\newcommand{\sym}{\mathcal{S}}

\newcommand{\qeds}{\qed\vspace{.2cm}}
\DeclareMathOperator{\tr}{Tr}

\DeclareMathOperator{\re}{Re}
\DeclareMathOperator{\imag}{Im}

\DeclareMathOperator{\im}{Im}
\DeclareMathOperator{\gram}{Gram}
\DeclareMathOperator{\cor}{\rm corank}

\DeclareMathOperator{\Ker}{Ker}

\title{
Universal completability, least eigenvalue frameworks, and vector colorings}

\author[1]{Chris Godsil}
\author[2]{David E.~Roberson}
\author[3]{Brendan Rooney}
\author[4]{Robert \v{S}\'{a}mal}
\author[5,6]{Antonios Varvitsiotis}
\affil[1]{Department of Combinatorics \& Optimization, University of Waterloo}
\affil[2]{Department of Computer Science, University College London}
\affil[3]{Department of Mathematical Sciences, KAIST}
\affil[4]{Computer Science Institute, Charles University}
\affil[5]{Centre for Quantum Technologies, National University of Singapore}
\affil[6]{School of Physical and Mathematical Sciences, Nanyang Technological University}

\date{}

\begin{document}
\maketitle
\begin{abstract}
\noindent An embedding  $i \mapsto p_i\in \R^d$ of  the vertices of a graph $G$ is called  \emph{universally completable} if the following holds: For any other embedding  $i\mapsto q_i~\in \R^{k}$ satisfying  $q_i\transpose q_j = p_i\transpose p_j$ for $i = j$ and $i$ adjacent to $j$, there exists an isometry mapping   the $q_i$'s to the $p_i$'s for all $ i\in V(G)$. The notion  of universal completability was  introduced recently  due to its relevance to the positive semidefinite matrix completion problem. In this work we focus on  graph embeddings  constructed using the eigenvectors  of the least eigenvalue of the adjacency matrix of $G$,  which we call    {\em least eigenvalue  frameworks}. We identify    two  necessary and sufficient conditions for such frameworks to be universally completable.  Our conditions also allow us to give   algorithms for determining whether a  least eigenvalue framework is universally completable. Furthermore, our  computations for Cayley graphs on  $\mathbb{Z}_2^n \ (n \le 5)$ show that almost all   of these graphs have universally completable least eigenvalue frameworks. 
In the second part of this work we study uniquely vector colorable (UVC) graphs, i.e., graphs for which the   semidefinite program corresponding to the Lov\'{a}sz theta number (of the complementary graph) admits a unique optimal solution. We  identify a sufficient condition for showing that a graph is UVC based on the universal completability of an associated framework. This allows us to prove that Kneser and $q$-Kneser graphs are UVC. Lastly, we show that  least eigenvalue frameworks of 1-walk-regular graphs always provide optimal vector colorings and furthermore,  we are able to 
 characterize all optimal vector colorings of such graphs. In particular, we give a necessary and sufficient condition for a 1-walk-regular graph to be uniquely vector~colorable. 
 \end{abstract}


\section{Introduction}\label{sec:intro}

%
%
%
%
%
%
%




A \emph{tensegrity graph} is  defined as a graph $G = ([n],E)$ where the edge set $E$ is partitioned  into   three disjoint sets $B$, $C$, and $S$. 
The elements of $B$, $C$, and $S$ are called  \emph{bars}, \emph{cables}, and \emph{struts} respectively. 
A \emph{tensegrity framework} $G(\p)$ consists of a tensegrity graph $G$ and an assignment of real vectors $\p = (p_1, p_2, \ldots, p_n)$ to the vertices of $G$.  In   this work we  consider   the vectors $p_i$ defining a tensegrity framework as   column vectors. For a tensegrity framework $G(\p)$ where   $p_i\in \R^{d}$ for all $i\in [n]$ we write $G(\p)\subseteq \R^d$. The {\em framework matrix} associated to a tensegrity framework $G(\p)\subseteq \R^d$, usually denoted by $P$,  is the $n\times d$ matrix whose $i^\text{th}$ row is given by the vector $p_i\transpose$. The {\em Gram matrix} of a tensegrity framework $G(\p)$, denoted by $\gram(p_1,\ldots,p_n),$  is defined as the  symmetric 
$n\times n$ matrix whose~$(i,j)$ entry is given by  $p_i\transpose p_j$ for all $i,j\in[n]$. Note  that for any tensegrity framework $G(\p)$ we have that   $\gram(p_1,\ldots,p_n)=PP\transpose$ and that ${\rm rank}(\gram(p_1,\ldots,p_n))=\dim \spn(p_1,\ldots,p_n).$

Consider two tensegrity frameworks $G(\p)$ and $G(\q)$. We say that   $G(\p)$  \emph{dominates}   $G(\q)$
if the following   three conditions hold:
\begin{itemize}
\item[$(i)$] $q_i\transpose q_j = p_i\transpose p_j,  \text{ for } ij \in B \text{ and for } i = j;  $
\item[$(ii)$] $q_i\transpose q_j \ge p_i\transpose p_j,  \text{ for } ij \in C;$
\item[$(iii)$] $q_i\transpose q_j \le p_i\transpose p_j,  \text{ for } ij \in S. $
\end{itemize}
We  use the notation    $G(\p)\succeq G(\q)$ to indicate  that  $G(\p)$ dominates $G(\q)$.
 Furthermore, two  tensegrity  frameworks $G(\p)$ and $G(\q)$ are called {\em congruent} if 
 \[\gram(p_1,\ldots,p_n)=\gram(q_1,\ldots,q_n).\]
 Lastly, 
  a tensegrity framework $G(\p)$ is called {\em universally completable}
if it is congruent to any framework it dominates, i.e., 
\[G(\p)\succeq G(\q) \ \Longrightarrow \   \gram(p_1,\ldots,p_n)=\gram(q_1,\ldots,q_n).\]

The notion of universal completability was originally introduced and investigated in~\cite{LV}. The main focus of this work is the study of least eigenvalue frameworks, i.e., graph embeddings obtained using the eigenvectors corresponding to the least eigenvalue of the adjacency matrix of the graph. Our first contribution is to identify two necessary and sufficient  conditions for showing that the least eigenvalue framework of a graph is   universally completable. The main tool in obtaining these characterizations is an explicit description  of the set of (Gram matrices) of tensegrity frameworks that are dominated by a framework $G(\p)$. The  proof of this fact is a slight modification of the proof of the main result from \cite{LV}.

Using these characterizations, we are able to identify a family of universally completable frameworks for the Kneser graph $K_{n:r}$ and the $q$-Kneser graph $qK_{n:r}$ for all $n \ge 2r+1$.
Furthermore, following other ideas from~\cite{LV}, we rephrase our second necessary and sufficient condition in terms of the Strong Arnold Property. This version of our result can be turned into an algorithm for deciding whether a least eigenvalue framework is universally completable. The algorithm only involves solving a system of $|V(G)|^2$ linear equations in $|E(\overline{G})|$ variables. Using this we show that 1293 out of the 1326 connected Cayley graphs for $\mathbb{Z}_2^n \  (n \le 5)$ have universally completable least eigenvalue frameworks.

In the second part of this paper we focus on {\em uniquely vector colorable} (UVC) graphs, i.e., graphs for which one of the semidefinite programs corresponding to the Lov\'asz theta number of the complementary graph admits a unique optimal solution. We give a sufficient condition for the optimality of a vector coloring $i\mapsto p_i$ by means of the universal completability of an appropriately defined tensegrity framework. Furthermore, we study the class of 1-walk-regular graphs for which we provide a full description of the set of optimal vector colorings. This yields a necessary and sufficient condition for a 1-walk-regular graph to be uniquely vector colorable.

The study of uniquely vector colorable graphs was initiated by Pak and Vilenchik in \cite{Pak} where the first family of UVC graphs was identified. The approach taken there is similar to ours, the main difference being that they use Connelly's sufficient condition for universal rigidity~\cite{C82}.  However, since vector colorings are defined in terms of inner products, the sufficient condition of Laurent and Varvitsiotis from~\cite{LV} can be applied in a much more direct and natural fashion. This is exemplified  by our ability to prove that a large and complex class of graphs, the Kneser and $q$-Kneser graphs, are uniquely vector colorable.

Some of the results in this paper were published as an extended abstract in the proceedings of the Eighth European Conference on Combinatorics, Graph Theory and Applications, EuroComb 2015~\cite{eurocomb}.

\subsection{Applications to other areas}

Laurent and Varvitsiotis originally introduced universal completability in~\cite{LV} due to its relevance to the low-rank positive semidefinite matrix completion problem. Recall that a (real) symmetric $n\times n$ matrix $Z$ is called {\em positive semidefinite} (psd), denoted by  $Z\in \sym^n_+$, if all of its eigenvalues are nonnegative.  To illustrate the connection, consider a tensegrity framework $G(\p)$,  let $X:=\gram(p_1,\ldots,p_n)$ be the corresponding Gram matrix and define 
\begin{align}\label{matrixconds}
\spec:=\Big\{ Z\in \sym^n_+:  Z_{ij} &= p_i\transpose p_j \text{ if } ij \in B \text{ or } i = j; \nonumber \\ 
Z_{ij} &\ge p_i\transpose p_j \text{ if } ij \in C;  \nonumber \\ 
Z_{ij} &\le p_i\transpose p_j \text{ if } ij \in S\Big\}. \nonumber
\end{align}
 The set  $\spec$  consists of all $n\times n$ psd  matrices   with the following properties: $(i)$ Diagonal entries and entries corresponding to bars coincide with the respective entries of $X$,    $(ii)$  entries corresponding to cables are  lower bounded by the respective entries of $X$, and  $(iii)$  entries corresponding to struts are upper bounded by the respective  entries   of $X$.  

As a matrix is psd if and only if it is the Gram matrix of a family of vectors, it follows that $\gram(p_1,\ldots,p_n)$ is an element of $\spec$ and that $G(\p)$ is universally completable if and only if this is the unique element~{of $\spec$}.  

For the remainder of this section  we consider  the special case where $S=C=~\emptyset.$ In this setting,   any tensegrity framework $G(\p)$ defines a  {\em $G$-partial matrix}, i.e., a matrix whose entries are only specified along  the diagonal and off-diagonal positions corresponding to edges of the $G$. In this case we have that 
\[\spec=\Big\{ Z\in \sym^n_+:  Z_{ij} = p_i\transpose p_j \text{ if } ij \in B \text{ or } i = j\Big\},\]
and   any element of   $\spec$ is  called a {\em psd completion} of the partial $G$-matrix specified  by the tensegrity framework $G(\p)$.  

Given a framework  $G(\p)$,  a question of fundamental interest   is to identify the smallest rank of an element in $\spec$. An  important instance of this problem  is the low-rank correlation matrix completion problem. The {\em correlation matrix} of a family of random variables $X_1, \ldots, X_n$ is the matrix whose $ij$-entry is equal to the correlation\footnote{The \emph{correlation} of two random variables is a rigorously defined measure of the degree to which those variables correlate with each other. The explicit definition is related to the covariance of the random variables and is outside the scope of this article. However, the details can be found in, e.g.~\cite{AThesis}.} of $X_i$ and $X_j$.  Equivalently, a matrix is a correlation matrix of some family of random variables if and only if it is psd with all diagonal entries equal to one~\cite{AThesis}. The rank of a correlation matrix turns out to be equal to the number of independent random variables. Consequently, identifying the smallest psd completion of a partial correlation matrix corresponds to finding the smallest dimensional model compatible with the observed correlations.

 The    {\em Gram dimension} of a graph $G$, denoted by ${\rm gd}(G),$   was introduced in~\cite{LV14}  to address the low-rank psd matrix completion problem described above. It  is defined as the smallest integer $k\ge 1$ with the following property: For any framework $G(\p)$ there exists an element $Z\in \spec$ such that ${\rm rank} (Z)\le k$. Notice    that  the Gram dimension  is  a well-defined graph parameter as it is bounded above by the number of vertices of $G$. Furthermore, it was shown in \cite{LV14} that ${\rm gd}(\cdot)$  is minor-monotone, i.e., whenever $H$ is a minor of $G$ we have that ${\rm gd}(H)\le {\rm gd}(G)$. By the graph minor theorem of Seymour and Robertson  it follows  that for any fixed integer $k\ge 1$, the class of graphs satisfying ${\rm gd}(G)\le k$ can be characterized by a finite list of   minimal forbidden minors. 
 The complete list of  forbidden minors was  for identified for   $k\in \{1,2,3,4\}$  in \cite{LV14}. 
 
  To show that a graph $G$ is a forbidden minor for  ${\rm gd}(H)\le k$ we need to construct  a framework $G(\p)$  for which any  $Z\in \spec$  satisfies ${\rm rank}(Z)>~k$.  Clearly,  placing a lower  bound on the rank of  all elements of $\spec$ is  a challenging task.  Nevertheless, there is one special case where this can achieved:  If  there exists a unique psd completion to the $G$-partial matrix specified by $G(\p)$, i.e.,  $\spec=\{\gram(p_1,\ldots,p_n)\}$. 
  Consequently, if $G(\p)$ is a universally completable framework with $\dim\spn(p_1,\ldots,p_n)>k$ then $G$ is a (not necessarily minimal) forbidden minor for the class of graphs satisfying 
${\rm gd}(G)\le k$. This is  the approach taken in \cite{LV}, where the notion of universal completability was introduced, to identify  families  of  forbidden minors for the Gram dimension.


The use of the term  ``universally completable" instead of, the perhaps more intuitive, ``uniquely completable" stems from a close relationship to theory of universal rigidity~\cite{C01}. 
A framework $G(\p)$ is called  \emph{universally rigid} if for any other framework $G(\q)\subseteq \R^d$ (for any $d\ge 1$)    
the following holds:
\[\|q_i - q_j\| = \|p_i - p_j\| \text{ for all } i \sim j \Longrightarrow \|q_i - q_j\| = \|p_i - p_j\| \text{ for all } i,j\in[n].\]
Here $\|\cdot\|$ denotes the usual Euclidean norm of a vector. 


Universal completability can be thought of as a ``spherical analog'' of universal rigidity, where norms of differences are replaced by inner products. In fact, the relationship can be made more precise, as was done in~\cite{LV}. Specifically, given a framework $G(\p)$, let  $G'$ by the graph obtained from $G$ by adding a vertex, labelled by 0, adjacent to every vertex in $G$, and let $\p' = (p_0, p_1, \ldots, p_n)$ where $p_0$ is the zero vector. Then it is not difficult to see that $G(\p)$ is universally completable if and only if $G'(\p')$ is universally rigid. One can also view the sufficient condition for universal completability given by Laurent and Varvitsiotis in~\cite{LV} as a spherical analog of the well known sufficient condition for universal rigidity due to Connelly~\cite{C82}.

\section{Universal completability}\label{sec:univcomp}
\subsection{Basic definitions and notation}\label{sec:definitions}

{\em Linear algebra.} We denote by $e_i$ the $i^\text{th}$ standard basis vector and by $\one$ the all-ones vector. Furthermore, we denote by $\spn(p_1,\ldots,p_n)$ the linear span of the vectors $\{p_i\}_{i=1}^n$.   The set of  $n\times n$ real symmetric matrices is denoted by $\sym^n$, and the set of matrices in $\sym^n$ with nonnegative eigenvalues, i.e.~the real positive semidefinite matrices, is denoted by $\sym^n_+$. Given a matrix $X\in \sym^n$ we denote its kernel/null space by $\Ker X$  and its image/column space by $\im X$. The  {\em corank} of a matrix $X\in \sym^n$, denoted $\cor X$, is defined as the dimension of~$\Ker X$. The {\em Schur} product of two matrices $X,Y\in \sym^n$, denoted by $X\circ Y$,  is the matrix whose entries are given by $(X\circ Y)_{ij}=X_{ij}Y_{ij}$ for all $i,j\in [n]$. Lastly, the {\em direct sum} of two matrices $X,Y\in \sym^n$, denoted by $X\oplus Y$,  is given by the  matrix
\[\begin{pmatrix} X& 0\\
 0 & Y\end{pmatrix}.\]
 Any  matrix $X\in \sym^n$ has real eigenvalues, and we denote the smallest one by $\lambda_{min}(X)$. 
 \medskip 
 
\noindent {\em Graph theory.} Given a graph $G=([n],E)$ we write $i\sim j$ to indicate that the vertices $i,j\in [n]$ are adjacent  and   $i\simeq j$ to indicate  that $i=j$ or $i\sim j$. For a vertex $i\in [n]$, we use $N[i]$ to denote the \emph{closed neighborhood} of $i$, i.e. $N[i]:=\{i\}\cup \{j \in [n] : j \sim i\}$. Also, for any subset $S\subseteq [n]$, we write $G\setminus S$ to denote the subgraph of $G$ induced by the vertices in $[n] \setminus S$.  We denote by $\overline{G}$ the {\em complement} of the graph $G$. A \emph{clique} in a graph is a subset of pairwise adjacent vertices, and an \emph{independent set} is a subset of pairwise nonadjacent vertices. A graph $G$ is called  \emph{split} if there exists a partition $\{C, I\}$ of the vertex set  such that $C$ is a clique in $G$ and $I$ is an independent set in $G$. The {\em Cayley graph} over a  group $\Gamma$ with inverse closed connection set $C\subseteq \Gamma \setminus \{\text{id}\}$ has the elements of $\Gamma$ as its vertices, such that vertices  $a, b\in \Gamma$ are adjacent if $a - b \in C$. If the group $\Gamma$ is abelian, then there exists a simple description of the eigenvalues and eigenvectors of the adjacency matrix of $G$. Specifically, if $\chi$ is a character of $\Gamma$, then the vector $(\chi(a))_{a\in \Gamma}$  is an eigenvector with corresponding eigenvalue $\sum_{c \in C}\chi(c)$. Moreover, this procedure provides  a full set of eigenvectors.

\subsection{A sufficient condition for universal completability}\label{sec:suffcond}

In this section we show that under appropriate assumptions  we can derive  a complete description for the set of tensegrity frameworks  that  a fixed tensegrity framework $G(\p)$ dominates (cf.~Theorem \ref{thm:tensegrity}). The proof of this fact is a slight generalization of the proof of the main theorem in \cite{LV}. As an immediate consequence, we identify a sufficient condition for showing that a tensegrity framework is universally completable as in \cite{LV}. Furthermore, this result is used in Section \ref{sec:vectcolor} where we study  uniquely vector colorable graphs.   

We start   with a definition which is central to the results in  this section.

\begin{defn}\label{def:sphstress}
Consider a tensegrity framework  $G(\p)\subseteq \R^d$  and let  $P
$ be the corresponding framework matrix. 
A \emph{spherical stress matrix} for $G(\p)$ is a  symmetric matrix  $Z \in \sym^n$ with the following properties:
\begin{enumerate}[(i)]
  \packed
  \item[$(i)$] $Z$ is positive semidefinite;
  \item[$(ii)$] $Z_{ij} = 0$ whenever $i \neq j, \ ij \not\in E$;
  \item[$(iii)$] $Z_{ij} \ge 0$ for all (struts) $ij \in S$ and $Z_{ij} \le 0$ for all (cables) $ij \in C$;
  \item[$(iv)$] $ZP = 0$;
  \item[$(v)$] $\cor(Z)=\dim \spn(p_1, \ldots, p_n)$.
  \end{enumerate}  
\end{defn} 

Although the above definition may make it seem like the existence of a spherical stress matrix is rare, the frameworks we will consider in this paper always admit a natural spherical stress matrix. Moreover, for any framework that is an optimal vector coloring (see Section~\ref{sec:vectcolor}), there always exists a (nonzero) matrix satisfying all but possibly condition $(v)$ of being a spherical stress matrix.

We continue with a simple lemma we use in the proof of our main result below. 

\begin{lemma}\label{lem:R}
Let $X\in \sym^n_+$ and $Y\in \sym^n$ satisfy $\Ker X\subseteq \Ker Y$. 
If $X=PP\transpose$ then there exists a symmetric matrix $R$ such that 
\[Y=PRP\transpose \text{ and } \ \im R\subseteq \im  P\transpose.\]
\end{lemma}
\proof We first prove the claim in the case that $P$ has full column rank. In this case we can extend $P$ to a full-rank square matrix $\tilde{P}$ and define the symmetric matrix  $\tilde{R}:=\tilde{P}^{-1}Y(\tilde{P}^{-1})\transpose$. By assumption we have that $\Ker X \subseteq \Ker Y$ and thus
\[\Ker X \oplus 0 = \Ker \tilde{P}(I\oplus 0)\tilde{P}\transpose \subseteq \Ker \tilde{P}\tilde{R}\tilde{P}\transpose = \Ker Y \oplus 0.\]
Since $\tilde{P}$ is invertible this gives  that $\Ker(I\oplus 0)\subseteq \Ker\tilde{R}$ and thus  $\tilde{R}=R \oplus 0.$
This shows that $Y=PRP\transpose$ for some symmetric matrix $R$. In this case we assumed that $P$ had full column-rank and so we have that $\im R \subseteq \im P\transpose$ since the latter is the whole space.


Lastly,  we consider the case when $P$ does not have full column rank. We have that $X = QQ\transpose$ for \emph{some} full column-rank matrix $Q$, and thus by the above there exists a symmetric matrix $R'$ such that $Y=QR'Q\transpose$.  Since $\im Q = \im X =\im P$ there exists a matrix $U$ such that $Q=PU$ and therefore $Y=PUR'U\transpose P\transpose$. If we let $E$ be the orthogonal projection onto $\im P\transpose$, then $EP\transpose = P\transpose$ and $PE = P$ since $E$ is symmetric. Thus $Y = PEUR'U\transpose EP\transpose$. Letting $R = EUR'U\transpose E$ completes the proof.\qeds

In~\cite{LV} it is shown that if there exists a spherical stress matrix for a framework $G(\p)$, where $\p$ spans the space it lives in, and $R= 0$ is the unique symmetric matrix satisfying condition $(b)$ in Theorem~\ref{thm:tensegrity} below, then $G(\p)$ is universally completable. Our main result of this section, presented here, slightly modifies the proof of this result to obtain a characterization of all frameworks dominated by $G(\p)$, assuming  that there exists a spherical stress matrix for $G(\p)$.


\begin{theorem}\label{thm:tensegrity}
Consider a tensegrity framework $G(\p)\subseteq \R^d$ and let  $P\in \R^{n\times d}$ be the corresponding framework matrix. Let $Z\in \sym_+^n$  be a spherical stress matrix for~$G(\p)$.  
The framework $G(\p)$ dominates the framework  $G(\q)$  
if and only if 
\begin{equation}\label{eq:gramequal}
{\rm Gram}(q_1, \ldots, q_n) = PP\transpose + PRP\transpose,
\end{equation}
where  $R$ is a symmetric $d\times d$ matrix satisfying: 
\begin{itemize}
\item[(a)] $\im R \subseteq \spn(p_1,\ldots,p_n)$;  
    \item[(b)] $p_i\transpose R p_j = 0 \text{ for } i = j \text{ and } ij \in B \cup \{\ell k\in C\cup S : Z_{\ell k} \ne 0 \};$
   \item[(c)]   $p_i\transpose R p_j \ge 0 \text{ for } ij \in C;$ 
     \item[(d)] $p_i\transpose R p_j \le 0 \text{ for } ij \in S.$
    \end{itemize}

\end{theorem}
\proof
Assume there exists a nonzero matrix $R\in \sym^d$  satisfying   $(a)$--$(d)$ and ${\rm Gram}(q_1, \ldots, q_n) = PP\transpose + PRP\transpose$. This shows that   $q_i\transpose q_j = p_i\transpose p_j + p_i\transpose R p_j$ for all $i,j\in [n]$ and using $(b)$--$(d)$  it follows that   $G(\p)\succeq~G(\q)$. For  the converse implication, say  that $G(\p)\succeq G(\q)$ and define  $X := PP\transpose = \text{Gram}(p_1, \ldots, p_n)$ and $Y := \text{Gram}(q_1, \ldots, q_n)$. Since $Z$ is a spherical stress matrix for $G(\p)$,  condition $(iv)$ implies  that $ZX = 0$. This shows  that $\im X \subseteq \Ker Z$. By $(v)$ we have  $\cor Z = {\rm rank}\  X$  and thus  $\Ker X = \im Z$. 
Furthermore, since  $G(\p)$ dominates $G(\q)$ and using the fact that $Y$ and $Z$ are positive semidefinite,  we have that
\begin{equation}\label{eqn}
0 \le \tr(ZY) = \sum_{ i \simeq j} Z_{ij}Y_{ij} \le \sum_{i \simeq j} Z_{ij}X_{ij} = \tr(ZX) = 0,
\end{equation}
and thus  \eqref{eqn} holds throughout with equality. In particular, again using that $Y$ and $Z$ are positive semidefinite, $\tr(YZ) = 0$ implies that  $YZ = 0$ and therefore $\Ker Y \supseteq \im Z = \Ker X$. If $v$ is a vector such that $Xv = 0$, then by the above $Yv = 0$ and therefore $(Y-X)v = 0$. This implies that $\Ker (Y - X) \supseteq \Ker X$. By Lemma~\ref{lem:R} we have that $Y-X = PRP\transpose$ for some symmetric matrix $R$ with $\im R \subseteq \im P\transpose = \spn(p_1, \ldots, p_n)$. By assumption $G(\p)\succeq G(\q)$  
 and thus $p_i\transpose R p_j \le 0 \text{ for } ij \in S$,  $p_i\transpose R p_j \ge 0 \text{ for } ij \in C$, and $p_i\transpose Rp_j = 0 \text{ for } ij \in B \text{ and } i =j$.  By \eqref{eqn} we have that $\sum_{i\simeq j}Z_{ij}(X_{ij}-Y_{ij})=0$ and since  $Z_{lk}(X_{lk}-Y_{lk})\ge 0$ for all $lk\in C\cup S$ it follows that 
  $X_{lk} = Y_{lk}$ for all $lk \in C \cup S$ with $Z_{lk}  \ne 0$.\qeds


\begin{remark}\label{rem:IRpsd}
Note that the matrix $R$ in the statement of Theorem~\ref{thm:tensegrity} satisfies $I+R \succeq 0$. This is because $\im R \subseteq \im P\transpose$ and $P(I+R)P\transpose$ is psd since it is a Gram matrix. Conversely, if $R$ is any matrix satisfying conditions $(a)$--$(d)$ and $I+R \succeq 0$, then $P(I+R)P\transpose$ is the Gram matrix of some set of vectors dominated by $\p$.
\end{remark}

As an immediate consequence of  Theorem \ref{thm:tensegrity}  we get the following sufficient condition for  showing  that a tensegrity framework $G(\p)$ is universally completable. This is essentially the same sufficient condition as that given in~\cite{LV}.



\begin{theorem}\label{thm:ucomp}
Consider  a tensegrity framework $G(\p)\subseteq \R^d$ and let $Z\in \sym^n_+$ be a  spherical stress matrix for $G(\p)$.  If $R=0$ is the unique symmetric matrix satisfying conditions $(a)$--$(d)$ of Theorem~\ref{thm:tensegrity}, then $G(\p)$ is universally completable. 
\end{theorem}

We note  that the existence of a spherical stress matrix is not a requirement for a framework to be universally completable. For an  example of a framework which is universally completable but admits no spherical stress matrix see~\cite{LV}.

\section{Least Eigenvalue Frameworks}\label{sec:lef}
\subsection{Definition and basic properties}\label{subsec:lef} 
In order to use  Theorem \ref{thm:ucomp}   to show that a  framework $G(\p)$  is universally completable, the first step is to construct a spherical stress matrix for $G(\p)$. In view of Definition \ref{def:sphstress},  it is not  obvious how to do this, and this task  is  equivalent to a feasibility semidefinite program with a rank constraint. In this section we show how for any graph $G$, using the eigenvectors of the  least eigenvalue of the adjacency matrix of $G$, we can construct a tensegrity framework which we call the {\em least eigenvalue framework}. Least eigenvalue frameworks are important  to this work as   they  come  with an associated spherical stress matrix. Consequently, to show that such a  framework is universally completable it suffices to check that the only matrix $R$ satisfying conditions $(a)$--$(d)$ of Theorem~\ref{thm:tensegrity} is $R = 0$.



\begin{defn}\label{def:lef}
Consider a  graph $G$ and let  $P$ be a matrix whose columns form an orthonormal basis for the eigenspace of the least eigenvalue of $G$. We say that the vectors $p_i$ that are the rows of $P$ are a \emph{least eigenvalue framework} of~$G$.
\end{defn}

Clearly,  there can be multiple least eigenvalue frameworks for a graph $G$, since there are many choices  of orthonormal basis for the least eigenspace. However, for any choice of orthonormal basis, the Gram matrix of the corresponding least eigenvalue framework will be equal to the orthogonal projector onto the least eigenspace of $G$. To see this, let $d$ be the dimension of the least eigenspace of $G$. Note that $P\transpose P = I_d$ and therefore $(PP\transpose)^2 = PP\transpose$. Therefore $PP\transpose$ is an orthogonal projector. Since the columns of $P$ are eigenvectors for the least eigenvalue of $G$, 
the column space of $PP\transpose$ is contained in the least eigenspace of $G$. To show that the column space of $PP\transpose$ is equal to the least eigenspace of $G$, it suffices to show that $PP\transpose$ has rank $d$. However, the eigenvalues of a projector take values in  $\{0,1\}$ and so its rank is equal to its trace. Therefore, ${\rm rank} (PP\transpose) = \tr(PP\transpose) =  \tr(P\transpose P) = \tr(I_d) = d$. Consequently, all least eigenvalue frameworks are congruent and thus indistinguishable for our purposes. We  refer to any such framework as \emph{the} least eigenvalue framework of~$G$.  

In general, one can consider eigenvalue frameworks for eigenvalues other than the minimum as well. In fact, this idea is not new, and such frameworks have been studied in their own right by algebraic graph theorists, under different names. Brouwer and Haemers~\cite{MR2882891} refer to eigenvalue frameworks as ``Euclidean representations," and use them to derive  the characterization of graphs with least eigenvalue at least $-2$, originally due to Cameron et al.~\cite{MR0441787}. Eigenvalue frameworks are used as a tool to prove statements about the structure of graphs from the geometry of its eigenspaces. In~\cite{MR1220704}, Godsil presents several results on distance-regular graphs that use this approach (Godsil refers to eigenvalue frameworks as ``representations"). Lastly, least eigenvalue frameworks also appear in the literature as the vertex sets of ``eigenpolytopes"~\cite{MR1638611}.

We will be interested in a similar, but more general, construction of frameworks based on eigenvectors  of a graph:
\begin{defn}
Consider a  graph $G$ and let $P$ be a matrix whose columns span the eigenspace of the least eigenvalue of $G$. We say that the vectors $p_i$ that are the rows of $P$ are a \emph{generalized least eigenvalue framework} of $G$.
\end{defn}

As an example of a generalized least eigenvalue framework, consider the projection $E_\tau$ onto the least eigenspace of some graph $G$. We claim that the rows, or columns of $E_\tau$ form a generalized least eigenvalue framework for $G$. Indeed, it is clear that these span the $\tau$-eigenspace of $G$, since this space is equal to the column space of $E_\tau$. Interestingly, the Gram matrix of this framework is $E_\tau^2 = E_\tau$, and so this framework is congruent to \emph{the} least eigenvalue framework of $G$, even though it was not explicitly constructed according to Definition~\ref{def:lef}. We will see another example of this phenomenon in Section~\ref{sec:kneser} when we construct a least eigenvalue framework for Kneser graphs.

In the next result   we show that for any  generalized  least eigenvalue framework there exists a canonical choice for a spherical stress matrix. 
\begin{proposition}\label{lem:stressmatrix}
Let $G$ be a tensegrity graph with no cables (i.e., $C=\emptyset$). Also let $A$ be the adjacency matrix of $G$ and set $\tau=\lambda_{min}(A)$. The matrix $A-\tau I$ is a spherical stress matrix for any generalized least eigenvalue framework of $G$. 
\end{proposition}

\proof We check the validity of conditions $(i)$--$(v)$ from Definition \ref{def:sphstress}.  Clearly, $A-\tau I \succeq 0$ and so condition $(i)$ holds. It is also trivial to see that conditions $(ii)$ and $(iii)$ hold. Condition $(iv)$ holds since the columns of $P$ are $\tau$-eigenvectors of $G$ by the definition of generalized least eigenvalue framework. For condition $(v)$, note that the corank of $A- \tau I$ is equal to the dimension of the $\tau$-eigenspace of $A$, which is exactly the dimension of the span of the  $\{p_i\}_{i=1}^n.$\qeds

Note that in  Proposition  \ref{lem:stressmatrix} we  restrict to tensegrities with no cables since the definition of the spherical stress matrix requires the entries corresponding to cables should  be less than 0. This   is clearly not satisfied for  $A-\tau I$. 

Given the restriction to tensegrity frameworks with no cables, the scope of   
 Proposition \ref{lem:stressmatrix}  appears to be limited. Nevertheless,  it turns out this is exactly what we need for the  study of uniquely vector colorable graphs in Section \ref{sec:vectcolor}.

\subsection{Conditions for universal  completability}
\label{sec:conditionsuc}
In this section we give a necessary and sufficient condition for showing that  a generalized least eigenvalue framework of a tensegrity graph with no cables is universally completable. We then apply our condition to show that the least eigenvalue framework of an odd cycle is universally completable.


\begin{theorem}\label{thm:lef}
Let $G$ be a tensegrity graph with $C=\emptyset$   and let $G(\p)\subseteq~\R^d$ 
be a generalized least eigenvalue framework of $G$. Then  $G(\p)$ is universally completable if and~{only if}
\begin{equation}\label{eq:conicinfty}
p_i\transpose R p_j = 0 \text{ for } i\simeq j 
\Longrightarrow R = 0,
\end{equation}
for all symmetric matrices $R\in \sym^d$ with $\im R \subseteq \spn(p_1, \ldots, p_n)$.
\end{theorem}
\proof
First, suppose $R\in \sym^d$ is a nonzero  matrix satisfying $p_i\transpose R p_j = 0$ for $i\simeq j$ and  $\im R \subseteq \spn(p_1, \ldots, p_n)$.
Without loss of generality we may assume that $\lambda_{min}(R)\ge -1$. Let $P$ be the framework matrix corresponding to $G(\p)$. Then the matrix $X:= P(I+R)P\transpose$ is positive semidefinite, and by assumption~satisfies
\[(P(I+R)P\transpose)_{ij} = p\transpose_ip_j + p_i\transpose R p_j = p\transpose_ip_j, \text{ for all } i\simeq j. \]
Furthermore, since $R \ne 0$ and $\im R \subseteq \spn(p_1, \ldots, p_n) = \im P\transpose$, the matrix $PRP\transpose$ is not zero. Thus  $X = P(I+R)P\transpose \ne PP\transpose$. Since $X$ is positive semidefinite, it is the Gram matrix of some set of vectors which form a framework that is dominated by $G(\p)$. Therefore, $G(\p)$ is not universally~completable.

For the other direction we use Theorem \ref{thm:ucomp} to show that $G(\p)$ is universally completable.  Let   $A$ be the adjacency matrix of $G$ and set  $\tau=\lambda_{min}(A)$. By Proposition~\ref{lem:stressmatrix} the matrix $A-\tau I $ is a spherical stress matrix for $G(\p)$.
Specializing the conditions $(a)$--$(d)$ from Theorem~\ref{thm:tensegrity} to $A-\tau I$ it remains to show that $R=0$ is the only symmetric matrix satisfying $p_i\transpose R p_j =~0$ for  $i\simeq j$ 
and $\im R \subseteq \spn(p_1, \ldots, p_n)$. This is exactly the assumption of the theorem.\qeds

%

\begin{remark}\label{rem:ucomplgenleast}
Notice that if $G(\p)$ is a  (generalized) least eigenvalue framework we have that 
\[p_i\transpose R p_j = 0 \text{ for } i\simeq j \Longleftrightarrow p_i\transpose R p_j = 0 \text{ for } i\sim  j.\]
To see this,  let   $A$ be the adjacency matrix of $G$, let $P$ be  the framework matrix associated to $G(\p)$ and $\tau=\lambda_{min}(A)$. By definition of the least eigenvalue framework we have that
$AP=\tau P$ which is equivalent to 
\[\tau p_i=\sum_{i\sim j}p_j, \text{ for all } i\in [n],\]
and the claim follows. 
\end{remark}

To illustrate the usefulness of Theorem \ref{thm:lef} we now show that the least eigenvalue framework of an odd cycle  is universally completable.
For  this we must first calculate  the eigenvectors corresponding to the least eigenvalue of an odd cycle. These are well known, but we briefly explain how to derive them.

The odd cycle of length $n:=2k+1$, denoted $C_{2k+1}$, can be described as the Cayley graph over  the abelian group $\mathbb{Z}_{n}:=\{0,\ldots,n-1\}$ with connection set $C = \{\pm 1 \mod n\}$.  As described in Section \ref{sec:definitions}    we see that $\lambda_{min}(C_{2k+1})=2\cos {2\pi k\over n } $ with multiplicity two. Furthermore,  a basis for the corresponding eigenspace is given by the vectors $\{ u,v\}\subseteq \mathbb{C}^n$ where $u(x)=\exp({2\pi i k x\over n})$ and $v(x)=\exp({2\pi i (n-k) x\over n})$. Recall that  the least eigenvalue framework of a graph is defined in terms of real vectors  but the  eigenvectors identified above are  complex valued.  It is easy to see that the vectors $\{ \re(u), \imag(u)\}\subseteq \R^n$ form a real valued orthogonal basis for  the least  eigenspace of $C_{2k+1}$. Note that $\re(u)_x=\cos({2\pi kx \over n})$ and $\imag(u)_x=\sin({2\pi kx \over n})$ for all $x\in \mathbb{Z}_n$. Consequently, the least eigenvalue framework for $C_{2k+1}$, up to a scalar,  is given by
\begin{equation}\label{eq:leasteigframcycle}
p_x=\left(\cos\Big({2\pi kx \over n}\Big),\sin\Big({2\pi kx \over n}\Big)\right)\transpose, \text{ for all } x\in \mathbb{Z}_n.
\end{equation}

%

%

Geometrically, this assigns $2k+1$ points evenly spaced around the unit circle to the vertices of $C_{2k+1}$ such that adjacent vertices are at maximum distance. 


\begin{theorem}\label{thm:ocuc}
The least eigenvalue framework of an odd cycle  is universally completable. 
\end{theorem}
\proof
Set $n:=2k+1$. Suppose $R \in \sym^2$ satisfies
\begin{equation}\label{eq:conic}
p_x\transpose R p_y = 0 \text{ for }x\simeq y.
\end{equation}
By Theorem \ref{thm:lef} it  suffices to show that $R = 0$.
For every $x \in  \mathbb{Z}_{n}$  it follows from  \eqref{eq:conic}   
 that $p_x$ is orthogonal to  the image of $\spn\{p_y : x \sim y\}$ under the map $R$. However, for   every $x\in \mathbb{Z}_{n}$ we have    $\spn\{p_y : x \sim y\}=~\R^2$ and thus  
 $p_x$ is orthogonal to  $\im R$. Since this is true for all $x \in \mathbb{Z}_n$, and  $\spn (p_x: x \in \mathbb{Z}_n)=\mathbb{R}^2$, we must have that $R = 0$.
\qeds


\subsection{Computations and the Strong Arnold Property}\label{sec:computations}

In the previous section we identified a necessary and sufficient condition for a  least eigenvalue framework to be universally completable (Theorem~\ref{thm:lef}). It is possible to turn this condition into an algorithm for determining when a least eigenvalue framework is universally completable. In fact, we investigate this approach in a follow-up to this paper~\cite{UVC2}. However, the algorithm investigated there only determines universal completability, it does not provide a method for determining all frameworks dominated by a least eigenvalue framework. Here we will present an alternative necessary and sufficient condition which presents a straightforward method for doing exactly this. The resulting algorithm corresponds to simply solving a homogeneous system of linear equations, and the framework is universally completable if and only if the system has only the trivial solution. Using this we  examine  how often the least eigenvalue framework of a graph happens to be  universally completable. Our computations  show that this is the case for the vast majority of  Cayley graphs on $\mathbb{Z}^2_n$ $(n\le 5)$.

To derive our second condition for the universal completability of least eigenvalue frameworks  we exploit a connection  with the Strong Arnold Property (SAP) that we now introduce. Consider  a graph $G = ([n], E)$ and let $A$ be its adjacency matrix. Set 
\[\mathcal{C}(G):= \{M \in \mathbb{R}^{n \times n} : M_{ij} = 0 \text{ if } i \not\simeq j\}.\]
A matrix $M \in \mathcal{C}(G)$ has  the \emph{Strong Arnold Property} if for every $X\in \sym^n$:
\[(A+I)\circ X=0
\text{ and }  MX = 0 \ \Longrightarrow \ X = 0.\]
The Strong Arnold Property is related to the celebrated  Colin de Verdi\`{e}re graph parameter \cite{cdv}, but Laurent and Varvitsiotis have also identified  a link between the SAP  and their sufficient condition for universal completability.

\begin{theorem}\cite{LV}\label{lem:LVSAP}
Consider a graph $G = ([n],E)$ and a matrix $M \in \mathcal{C}(G)$ with $\cor M=d$.  Let $P \in \mathbb{R}^{n \times d}$ be a matrix whose columns form an orthonormal basis for $\Ker M$ and let $p_1, \ldots, p_n$ denote the rows of $P$. The following are equivalent:
\begin{enumerate}
\item[(i)] $M$ has the Strong Arnold Property;
\item[(ii)] For any $d \times d$ symmetric matrix $R$,
\[p_i\transpose R p_j = 0 \text{ for all } i\simeq j  \Longrightarrow R = 0.\]
\end{enumerate}
\end{theorem}

We  now  give our second necessary and sufficient condition for universal completability. 


\begin{proposition}\label{lem:SAP}  
Let $G$ be a tensegrity graph with no cables, and let $G(\p)$ be its least eigenvalue framework. Furthermore, let $A$ be the adjacency matrix of $G$ and let $\tau=\lambda_{min}(A)$. The following are equivalent:
\begin{itemize}
\item[$(i)$] $G(\p)$ is universally completable;
\item[$(ii)$] $A-\tau I $ has the Strong Arnold Property. Explicitly, for any $X\in \sym^n$:
\[(A+I) \circ X = 0 \text{ and }  (A-\tau I)X = 0 \Longrightarrow X=0.\]
\end{itemize}
\end{proposition}

\proof The proof follows  by combining Theorem \ref{thm:lef} with Theorem \ref{lem:LVSAP}.\qeds



Proposition~\ref{lem:SAP} $(ii)$  provides us with a polynomial time algorithm for determining whether the least eigenvalue framework of a graph is universally completable, assuming we can compute its least eigenvalue exactly. In particular, finding all matrices $X$ such that $X_{ij} = 0$ when $i\simeq j$  and $(A - \tau I)X = 0$ is equivalent to solving a system of $|V(G)|^2$ equations (one for each entry of $(A-\tau I)X$) consisting of $|E(\overline{G})|$ variables. We then conclude that the least eigenvalue framework of $G$ is universally completable if and only if the only solution to this system of equations corresponds to $X = 0$.

Using Sage~\cite{sage}, we applied the above described algorithm\footnote{The Sage code will be made available to any interested party by email.} to all Cayley graphs for $\mathbb{Z}_2^n$ for $n \le 5$ and obtained the results summarized in the table below. Note that these graphs all have integral spectrum which allows us to use exact arithmetic when looking for possible solutions $X$.

\begin{table}[ht!]
\caption{Data for Cayley Graphs on $\mathbb{Z}_2^n$}
\centering
\[
\begin{array}{lrrrr}
n & \text{Num. Conn. Graphs} & \text{Num. Univ. Comp.}\\\hline
1 & 1 & 1\\
2 & 2 & 2\\
3 & 6 & 6\\
4 & 36 & 34\\
5 & 1326 & 1293
\end{array}
\]
\label{cubeliketable}
\end{table}

Note that we have not yet shown how to use this algorithm to determine all frameworks dominated by the least eigenvalue framework of a graph. To do this, we will need a correspondence between matrices $X$ satisfying the hypotheses of Proposition~\ref{lem:SAP} $(ii)$ and matrices $R$ satisfying $p_i\transpose R p_j = 0$ for $i \simeq j$, where $(p_1, \ldots, p_n)$ is the least eigenvalue framework of a graph. First we introduce some notation. For a graph $G$ with adjacency matrix $A$ and least eigenvalue $\tau$ with multiplicity $d$, we define the following:
\begin{align*}
\mathcal{X}(G) &= \{X \in \sym^n: (A+I)\circ X = 0 \text{ and } (A - \tau I)X = 0\} \\
\mathcal{R}(G) &= \{R \in \sym^d: p_i\transpose R p_j = 0 \text{ for } i \simeq j\}.
\end{align*}
Note that $\mathcal{X}(G)$ and $\mathcal{R}(G)$ are both clearly vector spaces. We will construct a linear map between these two spaces that will serve as the needed correspondence.

\begin{lemma}\label{lem:R2X}
Let $G$ be a graph and let $P$ be the framework matrix for the least eigenvalue framework of $G$. Define a map $\Phi$ as follows:
\[\Phi(R) = PRP\transpose \text{ for } R \in \mathcal{R}(G).\]
Then the map $\Phi$ is a linear bijection between $\mathcal{R}(G)$ and $\mathcal{X}(G)$, with $\Phi^{-1}(X) = P\transpose X P$. Moreover, $PP\transpose + \Phi(R)$ is psd if and only if $I+R$ is psd.
\end{lemma}
\proof
The fact that $\Phi$ is linear is obvious. Next we will show that $\Phi(R) \in \mathcal{X}(G)$ for all $R \in \mathcal{R}(G)$. If $R \in \mathcal{R}(G)$, then $R$ is symmetric and so $PRP\transpose$ is symmetric. Moreover,
\[\Phi(R)_{ij} = (PRP\transpose)_{ij} = p_i\transpose R p_j = 0 \text{ for } i \simeq j.\]
This implies that $(A+I) \circ \Phi(R) = 0$ where $A$ is the adjacency matrix of $G$. Also, by the definition of least eigenvalue framework, the columns of $P$ are eigenvectors for the minimum eigenvalue, say $\tau$, of $A$. Therefore,
\[(A-\tau I)\Phi(R) = (A- \tau I)PRP\transpose = 0,\]
as desired. This shows that the image of $\Phi$ is contained in $\mathcal{X}(G)$. Now we will show that $\Phi$ is surjective. Suppose $X \in \mathcal{X}(G)$. We have that 
\[\Phi(P\transpose X P) = PP\transpose X PP\transpose,\]
and moreover $PP\transpose =: E_\tau$ is the orthogonal projection onto the $\tau$-eigenspace of $A$. Since $(A- \tau I)X = 0$, the columns (and thus rows) of $X$ are all $\tau$-eigenvectors of $A$ and thus
\[\Phi(P\transpose X P) = E_\tau X E_\tau = X.\]
Therefore $\Phi$ is surjective.

Now we show that $\Phi$ is injective by verifying that $\Phi^{-1}(X) = P\transpose X P$. Since the columns $P$ are an orthonormal basis, we have that $P\transpose P = I$. Therefore,
\[P\transpose\Phi(R)P = P\transpose P R P\transpose P = R.\]
So $\Phi^{-1}(X) = P\transpose X P$.

Lastly, $PP\transpose + \Phi(R) = P(I+R)P\transpose$ which is psd if and only if $I+R$ is psd since $P$ has full column-rank.\qeds

The above lemma gives us the following corollary:

\begin{corollary}
Let $G$ be a tensegrity graph with no cables and let $\p$ be its least eigenvalue framework with corresponding framework matrix $P$. Then the matrix $PP\transpose + X$ is the Gram matrix of a framework dominated by $G(\p)$ if and only if $X \in \mathcal{X}(G)$ and $\lambda_{\min}(X) \ge -1$.
\end{corollary}
\proof
By Theorem~\ref{thm:tensegrity} and Remark~\ref{rem:IRpsd}, we have that $PP\transpose + X$ is the Gram matrix of a framework dominated by $G(\p)$ if and only if
\[PP\transpose + X = P(I+R)P\transpose\]
for some $R \in \mathcal{R}(G)$ such that $I+R$ is psd. By Lemma~\ref{lem:R2X} we have that this is equivalent to $X \in \mathcal{X}(G)$ and $PP\transpose + X \succeq 0$. The former implies that $\im X \subseteq \im PP\transpose$, and so the latter is equivalent to $\lambda_{\min}(X) \ge -1$ since $PP\transpose$ is an orthogonal projector.\qeds

By the above corollary, in order to determine all frameworks dominated by the least eigenvalue framework of a graph $G$, it suffices to determine $\mathcal{X}(G)$, which is equivalent to solving a system of linear equations. Then, for any $X \in \mathcal{X}(G)$, one can positively scale $X$ until $\lambda_{\min}(X) \ge -1$. This is especially useful when $\mathcal{X}(G)$ is 1-dimensional, which is a case we pay special attention to in~\cite{UVC2}.

%

\subsection{Two additional sufficient conditions}
\label{sec:othersuff}
In this section we use Proposition~\ref{lem:SAP}  to derive  two additional  sufficient conditions for the least eigenvalue framework of a graph to be universally completable.  For the  first of these, we recall  the following well-known result~(see, e.g., \cite{AGT}):

\begin{theorem}\label{lem:interlace}
Let $G$ be a graph and $H$ an induced subgraph of $G$. Then the least eigenvalue of $H$ is greater than or equal to the least eigenvalue of $G$.
\end{theorem}

Let $G$ be a tensegrity graph and let $\tau$ be  the minimum eigenvalue of its adjacency matrix. As an immediate consequence of Theorem \ref{lem:interlace} we have that 
 \begin{equation}\label{sferger}
 \lambda_{min}(G\setminus N[i])\ge \tau  \text{ for all } i\in V(G).
 \end{equation}
Using   \eqref{sferger} we next derive   a  sufficient condition for universal completability. 


\begin{proposition}\label{lem:nbrhdSAP}
Let $G=([n],E)$ be a tensegrity graph with no cables, let $A$ be its adjacency matrix and set  $\tau=\lambda_{min}(A)$.  
If 
\[\lambda_{min}(G\setminus N[i])>\tau\  \text{ for all } i\in V(G),\]
then the least eigenvalue framework  $G(\p)$ is universally completable.
\end{proposition}
\proof
 Suppose that $G(\p)$ is not universally completable.  We show there exists $i\in[n]$ for which $\lambda_{min}(G\setminus N[i]) = \tau$.  By Proposition~\ref{lem:SAP}, there exists a symmetric nonzero matrix $X\in \sym^n$ such that $(A+I) \circ X = 0$ and $(A-\tau I)X = 0$. Let $x_1, \ldots, x_n$ be the columns of $X$. Since $X$ is nonzero, there exists a vertex $i\in [n]$  such that $x_i\ne 0$. Furthermore, as $(A-\tau I)X = 0$ 
  it follows 
  that $x_i$ is a $\tau$-eigenvector of $G$. However, since $(A+ I) \circ X = 0$ we have that the entries of  $x_i$  corresponding to $N[i]$ are equal to zero.
  This implies that $G\setminus N[i]$ has $\tau$ as an eigenvalue.\qeds 


We continue  with our second sufficient condition for universal completability.

\begin{lemma}\label{lem:cliqueSAP}
Let $G=([n],E)$ be a tensegrity graph with no cables, let $A$ be its adjacency matrix and set  $\tau=\lambda_{min}(A)$.
Suppose there exists a clique $C$ in $G$ such  that the principal submatrix of $A-\tau I$ induced by the nodes in $[n] \setminus C$ is invertible. Then  the least eigenvalue framework $G(\p)$ is universally completable.
\end{lemma}
\proof
Suppose that $X\in \sym^n$ is a symmetric matrix satisfying $(A+I) \circ X = 0$ and $(A - \tau I)X = 0$. By Proposition \ref{lem:SAP} it suffices to show  that $X = 0$. By labeling the vertices of $G$ appropriately we can assume that $A-\tau I$ has the following block structure:
\[\begin{pmatrix}
D & B\\
B\transpose & F
\end{pmatrix},\]
where the upper left block corresponds to the clique $C$ whereas the lower right block corresponds to $[n]\setminus C$. 
By assumption we have that $F$ is invertible. Since $C$ is a clique, all of the entries of $A + I$ in the block corresponding to $C$ are 1. Since we require $(A+ I) \circ X = 0$, the corresponding block structure for $X$ is given by:
\[\begin{pmatrix}
0 & Z\\
Z\transpose & Y
\end{pmatrix}.\]
Then $(A-\tau I)X=0$ implies that 
\[FZ\transpose = 0 \text{ and } B\transpose Z + FY = 0.\]
Since $F$ is invertible it follows that $Z = Y = 0$ and thus $X = 0$.\qeds 

Using Lemma \ref{lem:cliqueSAP} we now show  that the least eigenvalue framework of any split graph is universally completable. 

\begin{corollary}
Let $G$ be a tensegrity graph with no cables. 
If $G$ is a split graph, then  the least eigenvalue framework $G(\p)$ is universally completable.
\end{corollary}
\proof
Since $G$ is a split graph there exists a partition $\{C,I\}$ of $V(G)$ such that $C$ is a clique and $I = V(G) \setminus C$ is an independent set of $G$. This implies that the principal submatrix of $A- \tau I$ corresponding to $V(G) \setminus C$ is a nonzero scalar multiple of the identity matrix, and is therefore invertible. By Lemma~\ref{lem:cliqueSAP} this implies that $G(\p)$ is universally completable.\qeds

\subsection{Application: Kneser Graphs}\label{sec:kneser}

In this section we  use Theorem \ref{thm:lef}  to show that a family of generalized least eigenvalue frameworks for the Kneser and $q$-Kneser graphs are universally completable. This is used later in Section \ref{sec:vectcolor} to show that for $n \ge 2r+1$,  both the Kneser graph $K_{n:r}$ and the $q$-Kneser graph $qK_{n:r}$ are uniquely vector colorable. Interestingly, even though the frameworks studied in this section are constructed as generalized least eigenvalue frameworks, we will later see that they are in fact congruent to the corresponding least eigenvalue framework.

 For two positive integers $n,r$ the {\em Kneser graph}, denoted by  $K_{n:r}$,  is the graph   whose vertices correspond to the subsets of $[n]$ of cardinality $r$,  
 where two vertices are adjacent  if the corresponding sets are disjoint. 
The {\em $q$-Kneser graph}, denoted by  $qK_{n:r}$,  has as its vertices the $r$-dimensional subspaces of the finite vector space $\mathbb{F}_q^n$, and two of these subspaces are adjacent if they are \emph{skew}, i.e., the subspaces intersect trivially.
In this section we construct  universally completable generalized least eigenvalue frameworks for  both $K_{n:r}$ and $qK_{n:r}$  (for $n \ge 2r+1$). We  only give the proof for the $q$-Kneser graphs, but the proof for  Kneser graphs is  similar; we  discuss the differences at the end of the section.


Let $P$ be a matrix with rows indexed by the $r$-dimensional subspaces of $\mathbb{F}_q^n$ (i.e.~by the vertices of $qK_{n:r}$) and columns indexed by the lines (1-dimensional subspaces) of $\mathbb{F}_q^n$ such that
\begin{equation}\label{eq:leasteigkneser}
P_{S,\ell}:= 
\begin{cases}
\alpha, & \text{ if } \ell \subseteq S, \\
\beta, & \text{ if } \ell \cap S = \{0\}.
\end{cases}
\end{equation}
In other words, $P$ is a  ``weighted incidence matrix'' of the $r$-dimensional subspaces and lines of $\mathbb{F}_q^n$. Further suppose that $\alpha$ and $\beta$ are chosen such that  $P\one=0$. 
The precise values of $\alpha$ and $\beta$ are not important (since we can apply a global scalar without changing the proof), but one suitable choice is
\[\alpha:= [r]_q - [n]_q\  \text{ and }  \
\beta := [r]_q,\]
where $[k]_q:= \frac{q^k - 1}{q-1} = \sum_{i=0}^{k-1} q^i$ gives   the number of lines contained in a $k$-dimensional subspace of $\mathbb{F}_q^n$. Using $P$ we construct a least eigenvalue framework for $qK_{n:r}$ by assigning to each $r$-dimensional subspace 
$S\subseteq \mathbb{F}_q^n$ the vector $p_S$ corresponding to the $S$-row of $P$. 
Note that the columns of $P$ are not orthogonal, however it is known that they span the least eigenspace of $qK_{n:r}$~\cite{newman06}. Therefore the vectors $\{p_S: S \in V(qK_{n:r})\}$ form  a generalized least eigenvalue framework of $qK_{n:r}$. Lastly,  note that the vectors $p_S$ lie  in $\mathbb{R}^{[n]_q}$ but do not span it, since they are all orthogonal to the all ones vector $\one$. Later we show  that $\spn(p_S: S \in V(qK_{n:r}))=\{\one\}^\bot$.

By Theorem~\ref{thm:lef}, 
to conclude  that $qK_{n:r}(\p)$ is universally completable it suffices to show that 
\[p_S\transpose Rp_T=0 \text{ for } S\simeq T \Longrightarrow R=0,\]
for any $R\in \sym^{[n]_q}$ with $\im R\subseteq \spn(p_S: S \in V(qK_{n:r}))$.
For  this we need to introduce some notation and an auxiliary lemma. 
For a subspace $F$ of $\mathbb{F}^n_q$, define $\one_F$ as a vector indexed by the lines of $\mathbb{F}_q^n$ as follows:

\begin{displaymath} 
\left(\one_F\right)_\ell = \left\{
\begin{array}{cl}
1, & \text{ if } \ell \subseteq F, \\
0, & \text{ if } \ell \cap F = \{0\}.
\end{array}
\right.
\end{displaymath}
We also define the following two subspaces of $\mathbb{R}^n$ for any subspace $F$ of $\mathbb{F}_q^n$:

\begin{equation}
\label{pande}
\begin{aligned}
P_F &:= \spn\left(\{p_S : S \cap F = \{0\}\} \cup \{\one\}\right), \text{ and }  \\
E_F &:= \spn\left(\{e_\ell : \ell \cap F = \{0\}\} \cup \{\one_F\}\right).
\end{aligned}
\end{equation}

We will need the following technical lemma:

\begin{lemma}\label{lem:subspaces}
Let $n \ge 2r+1$, and let $F$ be a subspace of $\mathbb{F}^n_q$ of dimension at most $r$. Then $P_F = E_F$.  
\end{lemma}
\proof
Clearly, $E_F$ is exactly the subspace of vectors which are constant on the lines contained in $F$. From this, it is easy to see that $P_F \subseteq E_F$. To show the other containment, we will prove that $e_\ell \in P_F$ for all $\ell$ skew to $F$, and that $\one_F \in P_F$.

First, suppose $\ell$ is a line skew to $F$. Then, since $n \ge 2r+1$, there exists some $r+1$ dimensional subspace $U$ of $\mathbb{F}^n_q$ containing $\ell$ and skew to $F$. Let $\mathcal{U}$ be the set of all $r$-dimensional subspaces of $U$. Since $U$ is skew to $F$, so is every element of $\mathcal{U}$. Therefore, for all $S \in \mathcal{U}$, we have $p_S \in P_F$. Furthermore, since $\one \in P_F$, we have that $\one_S \in P_F$ for all $S \in \mathcal{U}$. The vectors $\one_S$ for $S \in \mathcal{U}$ are all 0 on the lines not contained in $U$, and thus the matrix whose rows are these $\one_S$ vectors looks like
\[[M | \ 0 \ ],\]
where $M$ is the incidence matrix whose rows are indexed by the $r$-dimensional subspaces of $U$, and whose columns are indexed by the 1-dimensional subspaces of $U$. We will show that $e_{\ell'} \in P_F$ for all $\ell' \subseteq U$. To do this it suffices to show that $M$ has full column rank, which is equivalent to the matrix $M\transpose M$ having no zero eigenvalues.

To see this, note that $(M\transpose M)_{\ell'\ell''}$ is equal to the number of $r$-dimensional subspaces of $U$ which contain both $\ell'$ and $\ell''$. This value only depends on whether $\ell' = \ell''$ and is greater in the case where this holds. Therefore,
\[M\transpose M = aI + bJ\]
where $a > 0$ and $b \ge 0$. This clearly has no zero eigenvalues and thus $e_{\ell'} \in P_F$ for all $\ell' \subseteq U$, and in particular $e_\ell \in P_F$. Since $\ell$ was an arbitrary line skew to $F$, this shows that $e_{\ell'} \in P_F$ for all $\ell'$ skew to $F$.
To see that $\one_F \in P_F$, simply note that
\[\one_F = \one - \sum_{\ell \text{ skew to } F} e_\ell.\]\qeds

\begin{remark}
Setting $F$ to be   the zero subspace, Lemma \ref{lem:subspaces} implies  that
\[\spn(\{p_S : S \in V(qK_{n:r})\} \cup \{\one\})=\mathbb{R}^{[n]_q},\]
and since $p_S\transpose \one=0 $ for all $S\subseteq \mathbb{F}_q^n$ it follows that 
 \begin{equation}\label{eq:span}
 \spn(p_S : S \in V(qK_{n:r}))= \{\one\}^\perp.
 \end{equation}
This is used  in the proof of our main theorem below.
 \end{remark} 

\begin{theorem}\label{thm:Kneser}
For $n \ge 2r+1$, the generalized least eigenvalue framework of $qK_{n:r}$ described in  \eqref{eq:leasteigkneser}   is universally completable.
\end{theorem}
\proof Suppose that $R$ is an $[n]_q \times [n]_q$ symmetric matrix satisfying $\im R \subseteq \{\one\}^\bot$ and
\begin{equation}\label{defergrh}
p_S\transpose R p_T = 0 \text{ for all } S,T \in V(qK_{n:r}) \text{ such that } S \sim T.
\end{equation}
Combining Theorem \ref{thm:lef} with \eqref{eq:span}, it  suffices to show that $R=0$. For any $T \in V(qK_{n:r})$ it follows from \eqref{defergrh} that  the vector $R p_T$  is orthogonal to $p_S$ for all $S \in V(qK_{n:r})$ skew to $T$. Furthermore, as $\im R \subseteq \{\one\}^\bot$ the vector $Rp_T$ is orthogonal to $\one$. This implies that $Rp_T$ is orthogonal to $P_T$ as defined in \eqref{pande}. By Lemma~\ref{lem:subspaces}, we have that $P_T = E_T$ and therefore
\[Rp_T \perp e_\ell \text{ for all } \ell \text{ skew to } T.\]
Since $R$ is symmetric, this implies that
\[Re_\ell \perp p_T \text{ for all } T \text{ skew to } \ell.\]
As $\im R \subseteq \spn(\one)^\bot$, the latter implies that $Re_\ell$ is orthogonal to $P_F$ for $F = \ell$. Applying Lemma~\ref{lem:subspaces} again, we obtain
\[Re_\ell \perp E_\ell = \mathbb{R}^{[n]_q}.\]
Since this is true for all lines $\ell$ of $\mathbb{F}_q^n$, we have that $R = 0$.\qeds

The proof for Kneser graphs (both the lemma and the theorem) is essentially identical to the above, except that subspaces are replaced by subsets and lines are replaced by elements. Therefore we have the following:
\begin{theorem}
For $n \ge 2r+1$, the generalized least eigenvalue framework of $K_{n:r}$ described in  \eqref{eq:leasteigkneser} is universally completable.
\end{theorem}


\section{Vector Colorings}\label{sec:vectcolor}
\subsection{Definitions and properties}

For $t \ge 2$, a {\em vector $t$-coloring} of a graph $G=([n],E)$  consists of  an assignment $\p = (p_1, \ldots, p_n)$ of real {\em unit} vectors to the vertices of $G$ such that  
\begin{equation}\label{eq:vcoloring}
p_i\transpose p_j \le \frac{-1}{t-1},  \text{ for all  } i \sim j.
\end{equation}
 We say that $\p$ is a {\em strict vector $t$-coloring} if \eqref{eq:vcoloring}   holds with equality for all edges of $G$. 
The {\em vector chromatic number} of a nonempty graph $G$, denoted $\chiv(G)$, is the minimum real number $t\ge 2$ such that $G$ admits a vector $t$-coloring. The vector chromatic number of empty graphs is defined to be one. The {\em strict vector chromatic number}, $\chisv(G)$, is defined analogously. We say that the \emph{value} of a vector coloring is the smallest $t$ for which~(\ref{eq:vcoloring}) is satisfied.

Vector and strict vector colorings, as well as their associated chromatic numbers, were introduced  by Karger, Motwani, and Sudan in~\cite{KMS}. As it turns out, for any graph $G$ we have  $\chisv(G)=\vartheta(\overline{G})$, where $\vartheta(\cdot)$ denotes the Lov\'{a}sz theta number of a graph \cite{lovasz}. Furthermore,  for any graph $G$  we have that
$\chiv(G)=\vartheta'(\overline{G})$, where $\vartheta'(\cdot)$ is a variant of the Lov\'{a}sz theta function introduced by Schrijver in~\cite{schrijver}.  

 Clearly, for any graph $G$ we have that $\chiv(G)\le \chisv(G)$.  Furthermore, notice 
 that if $G$ admits  a $k$-coloring (in the usual sense), then mapping each color class to a vertex of the $k-1$ simplex centered at the origin corresponds to a  strict vector $k$-coloring. This implies that $\chisv(G) \le \chi(G)$.
 
Our main goal in this section is to  identify families of graphs  that admit  unique vector (resp. strict) vector colorings.  
Similarly to    tensegrity frameworks, any (strict) vector $t$-coloring can be rotated, reflected, or otherwise orthogonally transformed to generate  another (strict) vector $t$-coloring. This is analogous to permuting the colors in a usual coloring of a graph.  Consequently,  when 
 defining uniquely vector colorable graphs we must  mod out by this equivalence  to arrive at  a meaningful definition. As with tensegrity frameworks, we accomplish this using Gram matrices.

\begin{defn}
A graph $G$ is called \emph{uniquely (strict) vector colorable} if for any two \emph{optimal} (strict) vector colorings $\p$ and $\q$, we have 
\[\gram(p_1, \ldots, p_n) = \gram(q_1, \ldots, q_n).\]
\end{defn}

In this section we identify an interesting  connection between universal completability and uniquely vector colorable graphs.  Specifically, in Section \ref{sec:unicity} we develop a sufficient condition for showing that a vector coloring of a graph $G$ is the unique optimal vector coloring of $G$
Combining this with our results from  the previous section allows us to  show that for $n \ge 2r+1$,  both the Kneser graph $K_{n:r}$ and the $q$-Kneser graph $qK_{n:r}$ are uniquely vector colorable. Furthermore, in Section~\ref{sec:1walkreg} we introduce the class of 1-walk-regular graphs for which we can fully characterize the set of optimal vector colorings.  To achieve this we  need to use   the characterization of frameworks dominated by a fixed tensegrity framework $G(\p)$ given in Theorem \ref{thm:tensegrity}.

\subsection{Uniqueness of vector colorings}\label{sec:unicity}

In this section  we  mainly focus  on unique vector colorability, but for the graph classes we consider this is    equivalent to unique strict vector colorability.  In order to show  that a graph $G=([n],E)$ is uniquely vector colorable we   start with a candidate  vector coloring $\p$ of $G$ and show it is the unique optimal vector  coloring of $G$.
To achieve  this   we use the tools we  developed in the previous sections  concerning tensegrity frameworks.
Specifically,  we associate to the  vector coloring  $\p$   a  tensegrity framework $\tilde{G}(\p)$ and relate the universal completability  of $\tilde{G}(\p)$ to the uniqueness of  $\p$ as an optimal  vector coloring. 

\begin{defn}\label{def:tens}
Consider a graph $G=([n],E)$ and let $\p$ be a  vector coloring of $G$. Define $\tilde{G}$ to be  the tensegrity graph obtained by $G$  by setting $S=E$ (and thus $C=B=\emptyset$)  and let $\tilde{G}(\p)$   the corresponding tensegrity~framework.
\end{defn}

Using the construction  described above we now state and prove a lemma that is  used   throughout this section.  
\begin{lemma}\label{lem:bettervalue}
Consider a graph $G=([n], E)$ and let  $\p$  be a \emph{strict} vector coloring  of  $G$. The vector coloring    $\q$ achieves better (i.e., smaller) or equal value compared to the vector coloring $\p$ if and only if   
 $\tilde{G}(\p)\succeq \tilde{G}(\q)$.
\end{lemma} 

\proof Let $t$ be the value of  $\p$ as a strict vector coloring. 
Let  $\q$  be a vector coloring of $G$  satisfying $q_i\transpose q_j\le {-1\over t'-1}$ for some $t'\le t$. As $\p$ is a strict vector $t$-coloring we have  that $q_i\transpose q_j\le{-1\over t-1}= p_i\transpose p_j$ for all $i\sim j$. Since the vectors $\{p_i\}_{i=1}^n$ and $\{q_i\}_{i=1}^n$ have unit norm it follows that $\tilde{G}(\p)\succeq \tilde{G}(\q)$.  The converse direction follows easily (and holds even if $\p$  is not strict as a vector coloring).\qeds


\begin{remark}
Note that the forward implication in  Lemma \ref{lem:bettervalue} is not guaranteed to hold if $\p$ is not a strict vector $t$-coloring.  Indeed, in this case there exist adjacent vertices $i$ and $j$ such that $p_i\transpose p_j < -1/(t-1)$ so there might exist a vector coloring $\q$ whose value is better compared to $\p$ but satisfies
\[p_i\transpose p_j < q_i\transpose q_j < \frac{-1}{t-1}.\]
This  shows  that $\tilde{G}(\p)\not \succeq \tilde{G}(\q)$.
\end{remark}


Using Lemma~\ref{lem:bettervalue} we arrive at our main result in this section where we make the connection between universal completability and unique vector colorability.

\begin{theorem}\label{thm:conditions}
Consider a graph $G=([n], E)$, let $\p$ be a strict  vector coloring and let $\tilde{G}(\p)$ be the tensegrity framework given in Definition~\ref{def:tens}. Then:
\begin{itemize}
\item[(i)] If $\tilde{G}(\p)$ is universally completable, then $\p$ is the unique optimal vector coloring of $G$.
\item[(ii)] If $\p$ is optimal as a vector coloring, then $G$ is uniquely vector colorable if and only if $\tilde{G}(\p)$ is universally completable.
\end{itemize}
\end{theorem}

\proof 
\textit{(i)}  
Let $\q$ be a vector coloring of $G$  whose value is better or equal compared to $\p$.  By Lemma~\ref{lem:bettervalue} we have that $\tilde{G}(\p)\succeq \tilde{G}(\q)$. By assumption   $\tilde{G}(\p)$ is universally completable which implies  $\gram(p_1,\ldots,p_n)=\gram(q_1,\ldots,q_n)$.

\textit{(ii)} 
By assumption we have  $p_i\transpose p_j=-1/(\chiv(G)-1)$  for all $i\sim j$.  First, assume  that $G$ is uniquely vector colorable and consider a framework $\tilde{G}(\q)$ satisfying  $\tilde{G}(\p)\succeq\tilde{G}(\q) $.   This gives   that   $q_i\transpose q_j\le p_i\transpose p_j= -1/(\chiv(G)-1)$ for all $i\sim j$ and  since $\p$ is the unique optimal vector coloring of $G$ it follows that $\gram(p_1,\ldots,p_n)=\gram(q_1,\ldots,q_n)$. Conversely, say that  $\tilde{G}(\p)$  is universally completable and let $\q$ be another optimal vector coloring of $G$.  This means that $q_i\transpose q_j\le -1/(\chiv(G)-1)=p_i\transpose p_j$ for all $i \sim j$. Since $\tilde{G}(\p)$ is universally completable it follows that $\gram(p_1,\ldots,p_n)=\gram(q_1,\ldots,q_n)$.\qeds 

We conclude this section with  an easy  application of Theorem \ref{thm:conditions}  where we show that odd cycles, Kneser, and $q$-Kneser graphs are uniquely vector colorable. 

\begin{theorem}
For $k \in \mathbb{N}$, the odd cycle $C_{2k+1}$ is uniquely vector colorable.
\end{theorem}
\proof
Let $\p$ be the least eigenvalue framework of $C_{2k+1}$, as described in Section~\ref{sec:conditionsuc}. It is easy to see that the inner products of the vectors in $\p$ are constant on edges and that this constant is negative. Therefore, after appropriate scaling, $\p$ is a strict vector coloring. By Theorem~\ref{thm:conditions}, it remains to show that $\tilde{C}_{2k+1}(\p)$ is universally completable. However this was shown in Theorem~\ref{thm:ocuc}.\qeds

\begin{theorem}
For $n \ge 2r+1$,  both the Kneser graph $K_{n:r}$ and the $q$-Kneser graph $qK_{n:r}$ are uniquely vector colorable.
\end{theorem}
\proof
Let $\p$ denote the generalized least eigenvalue framework for  the $q$-Kneser graph $qK_{n:r}$ described in Section \ref{sec:kneser}. It is easy to see that this framework has constant inner product on edges of $qK_{n:r}$, and a straightforward computation shows that this constant is negative. Therefore, after appropriate scaling, this forms a strict vector coloring. Consequently, by Theorem~\ref{thm:conditions} $(i)$, it suffices to show that $q\tilde{K}_{n:r}(\p)$ is universally completable.  That was already established  in Theorem~\ref{thm:Kneser}. The case of Kneser graphs follows in a similar manner.\qeds

For $n \le 2r-1$, the graphs $K_{n:r}$ and $qK_{n:r}$ are empty, and so they are clearly not uniquely vector colorable. Furthermore, for $n = 2r$ (and $r > 1$), the graph $K_{n:r}$ is disconnected and therefore not uniquely vector colorable. This leaves the case $n=2r$ for the $q$-Kneser graphs. These graphs are not bipartite, so it is not clear if they are uniquely vector colorable. However, using the algorithm described in Section~\ref{sec:computations}, we found that $2K_{4:2}$ is not uniquely vector colorable.



\subsection{1-Walk-Regular Graphs}\label{sec:1walkreg}

In this section we focus on  the class of 1-walk-regular graphs. These graphs  are relevant to this work since they exhibit  sufficient regularity so as to guarantee that their least eigenvalue frameworks are always (up to a global scalar)   strict vector colorings.  This fact  combined with  Theorem \ref{thm:tensegrity} allows us to characterize the set of optimal vector colorings for a  1-walk-regular graph (cf. Theorem \ref{cor:1walkreg}). 
This implies   that   the least eigenvalue framework of a 1-walk-regular graph  is always an optimal strict vector coloring and moreover, yields  a necessary and sufficient condition for a 1-walk-regular graph to be uniquely vector colorable.

\begin{defn}  
A graph with adjacency matrix $A$ is called  {\em 1-walk-regular} if there exist $a_k, b_k \in \mathbb{N}$ for all $k \in \mathbb{N}$ such that:
\begin{enumerate}
\item[$(i)$] $A^k \circ I = a_k I$;
\item[$(ii)$] $A^k \circ A = b_k A$.
\end{enumerate}
\end{defn}
Equivalently, a graph is \emph{1-walk-regular} if for all $k \in \mathbb{N}$, $(i)$ the number of walks of length $k$ starting and ending at a vertex does not depend on the choice of vertex, and $(ii)$ the number of walks of length $k$ between the endpoints of an edge does not depend on the edge.

Note that a 1-walk-regular graph must be regular. Also, any graph which is vertex- and edge-transitive is easily seen to be 1-walk-regular. Other classes of 1-walk-regular graphs include distance regular graphs and, more generally, graphs which are a single class in an association scheme.


Our main result in this section is a characterization of the set of optimal vector colorings of a 1-walk-regular graph which we now state and prove. 
\begin{theorem}\label{cor:1walkreg}
Consider  a 1-walk-regular graph $G=([n],E)$. Let $G(\p)\subseteq \R^d$ be its least eigenvalue framework  and 
 $P\in \R^{n\times d}$  the corresponding framework matrix. 
 The  vector coloring $\q$   is optimal  if and only if
\begin{equation}\label{eq:optimalveccoloring1walk}
\gram(q_1, \ldots, q_n) = \frac{n}{d}\left(PP\transpose + PRP\transpose\right),
\end{equation}
where $R\in \sym^d $ is a symmetric matrix satisfying
\[p_i\transpose R p_j = 0 \text{ for all  } i \simeq j.\]
\end{theorem}

\proof
Let $A$ be the adjacency matrix of $G$ and set  $\tau:=\lambda_{\min}(A)$. Note that $d = \cor(A-\tau I)$, and  that  the matrix   $E_\tau:= PP\transpose$  is the orthogonal projector  onto the $\tau$-eigenspace of $G$. We  first show    that the assignment   $i\mapsto  \sqrt{\frac{n}{d}}p_i$ is   a strict vector coloring. For this  we show that  $p_i\transpose p_i = d/n$ for all $i \in [n]$ and that $p_i\transpose p_j$ is a negative constant for all $i \sim j$.

First, note that  $E_\tau$ can be expressed  as a polynomial in $A$. To see this set  
\[Z:=\prod_{\lambda \ne \tau}\frac{1}{\tau - \lambda}(A - \lambda I),\]
where the product is over eigenvalues of $A$. Considering $Zv$ where $v$ ranges over an orthonormal basis composed of eigenvectors of $A$ shows that $Z = E_\tau$. 

Since $G$ is 1-walk-regular and $E_\tau$ is a polynomial in $A$  there  exist scalars  $a, b$ such that 
\begin{equation}\label{eq:alconditions}
E_\tau \circ I = aI \ \text{ and }  \ E_\tau \circ A = bA.
\end{equation}\
Using the fact that  $E_\tau=PP\transpose$ it follows from \eqref{eq:alconditions} that 
\begin{equation}
p_i\transpose p_i=a, \text{ for all } i\in [n], \text{ and } p_i\transpose p_j=b \text{ for all } i\sim j.
\end{equation}
Moreover, since $E_\tau$ is the  projector onto $\Ker(A-\tau I)$ and $d=\cor(A-\tau I)$, we have that $\tr(E_\tau)={\rm rank}(E_\tau)=~d $. On the other hand Equation  \eqref{eq:alconditions} implies that $\tr(E_\tau)=na$ and thus  $a=d/n$, as previously claimed.

Let $\text{sum}(M)$ denote the sum of the entries of the matrix $M$. Using \eqref{eq:alconditions} combined with the fact that  $G$ is $r$-regular for some $r$, we get that 
\[
brn = \text{ sum}(A \circ E_\tau) = \tr(AE_\tau) = \tr(\tau E_\tau) = \tau d,  
\] 
and thus $b=\tau d/nr < 0$, since $\tau < 0$. 

For $i\in [n]$, set $ \tilde{p}_i:= \sqrt{\frac{n}{d}}p_i$. Since 
\[\tilde{p_i}\transpose \tilde{p}_i=1,  \text{ for all } i\in [n]   \text{ and } \ \tilde{p_i}\transpose \tilde{p}_i={\tau\over r} \text{ for all } i\sim j,\]
the assignment $\tilde{\bf p}$ 
is a strict  vector coloring of $G$.

Let  $\tilde{G}$ denote the tensegrity graph obtained from the graph $G$  by making all its  edges  into struts (recall  Definition \ref{def:tens}).  
By Lemma \ref{lem:bettervalue} we have that  a vector coloring $\q$  of  $G$    achieves the  same or better   value compared   to $\tilde{\bf p}$ if and only if  
 $\tilde{G}(\tilde{\bf p})\succeq \tilde{G}(\q)$. 
Furthermore, since $\p$ was the least eigenvalue framework of $G$, by Proposition~\ref{thm:lef} we know that  $A-\tau I$ is a spherical stress matrix for $\tilde{G}(\tilde{\p})$. Consequently, it follows from   Theorem~\ref{thm:tensegrity}  that  $\tilde{G}(\tilde{\bf p})\succeq \tilde{G}(\q)$ is equivalent to 
\begin{equation}\label{eq:allvectorcoloring}
\gram(q_1,\ldots,q_n)= \frac{n}{d}\left(PP\transpose + PRP\transpose\right),
\end{equation}
for some  $R\in \sym^d$  satisfying    $  p_i\transpose Rp_j = 0$ whenever $i \simeq  j$. 
Lastly, this implies that  
\[q_i\transpose q_j=\tilde{p}_i\transpose\tilde{p}_j \text{ for all } i\simeq j,\]
and thus $\q$ is a vector coloring of $G$ achieving the same value as the vector coloring $\tilde{\bf p}$.
\qeds

As a  consequence of Theorem \ref{cor:1walkreg}   we identify  a necessary and sufficient     condition for showing that 
a 1-walk-regular graph is uniquely vector~colorable.

\begin{corollary}\label{thm:1hom}
Let $G=([n],E)$ be 1-walk-regular and let $G(\p)\subseteq \R^d$ be its least eigenvalue framework. We have that:
\begin{itemize}
\item[$(i)$] The  assignment  $i\mapsto \sqrt{\frac{n}{d}}p_i $ is an optimal strict vector coloring of $G$.
\item[$(ii)$]    $G$ is uniquely vector colorable if and only if for any     $R\in \sym^d $ we have:
\[p_i\transpose R p_j = 0 \text{ for all  } i \simeq j \Longrightarrow R=0.\]
\end{itemize}
\end{corollary}

We note that the  construction of an optimal strict vector coloring described above   originally appeared  in~\cite{sabvshed}, though the proof of optimality  there is different and uniqueness is not discussed. Since the Kneser and $q$-Kneser graphs are 1-walk-regular and UVC, their unique vector colorings must be their least eigenvalue frameworks (up to a scalar multiple). This means that the frameworks constructed for these graphs in Section~\ref{sec:kneser} must have been congruent to their least eigenvalue frameworks, even though they only appeared to be generalized least eigenvalue frameworks by description.






\section{Concluding remarks}\label{sec:discuss}

In the first part of this work we considered general tensegrity frameworks. We showed that for any framework $G(\p)$ for which there exists a spherical stress matrix,  we can provide a description of the set of frameworks that are dominated by $G(\p)$. We then introduced least eigenvalue frameworks and identified two necessary and sufficient conditions for determining when a least eigenvalue framework is universally completable. Using these conditions we showed that a family of least eigenvalue frameworks for the Kneser and $q$-Kneser graphs are universally completable. Lastly, by reformulating our conditions in terms of the Strong Arnold Property, we gave an efficient algorithm for determining when a least eigenvalue framework is universally completable. Using this, we collected data on Cayley graphs over  $\mathbb{Z}_2^n\  (n\le 5)$ indicating that it is fairly common for the least eigenvalue framework to be universally completable.


In the second part of this work, we introduced the notion of unique vector colorability and showed that certifying the optimality of a vector coloring that is strict can be reduced to the universal completability of an appropriately defined tensegrity framework. This fact allowed us to conclude that odd cycles, Kneser and $q$-Kneser graphs are uniquely vector colorable. Lastly, we characterized the set of optimal vector colorings for the class of 1-walk-regular graphs.   As a corollary we got that the least eigenvalue framework of a 1-walk-regular graph is always an optimal strict vector coloring and moreover, we identified a necessary and sufficient condition for a 1-walk-regular graph to be uniquely vector colorable.


To our knowledge, the only other class of graphs known to be uniquely vector colorable is due to Pak and Vilenchik~\cite{Pak}. There, they identify a sufficient condition (on the eigenvalues of a regular graph $G$) for the categorical product of $G$ and a complete graph to be uniquely vector colorable. Furthermore, they  give an explicit construction of a class of uniquely vector colorable graphs.

As a follow up to this work, we are currently preparing two other articles on vector colorings and unique vector colorability~\cite{UVC2, UVC3}. One of these focuses on the relationship between vector colorings and graph homomorphisms, and in particular, how to use knowledge of the former to obtain information about the latter. The other article investigates vector colorings of categorical products. Here we prove a vector coloring analog of the well known Hedetniemi conjecture, and also greatly generalize the result of Pak and Vilenchik mentioned above.

\paragraph{Acknowledgements:}  D.~E.~Roberson and A.~Varvitsiotis were supported in part by the Singapore National Research Foundation under NRF RF Award No. NRF-NRFF2013-13. R.~\v{S}\'{a}mal was partially supported by grant GA \v{C}R P202-12-G061 and by grant LL1201 ERC CZ of the Czech Ministry of Education, Youth and Sports.

\bibliographystyle{plainurl}

\bibliography{main.bib}

\end{document}